\documentstyle[amscd,amssymb,12pt]{amsart}
\input xypic \xyoption{curve}
\input epsf

\def\hcorrection#1{\advance\hoffset by #1 }
\def\vcorrection#1{\advance\voffset by #1 }

\vcorrection{-.7in}
\hcorrection{-.5in}

\newcommand{\B}[1]{{\bold#1}} 
\newcommand{\C}[1]{{\cal#1}} 

\theoremstyle{plain}

\theoremstyle{definition}

\theoremstyle{definition}

\theoremstyle{remark}

\numberwithin{equation}{section}

\begin{document}

\pagestyle{plain}
\addtolength{\footskip}{.3in}

\title{Remarks on\\
Quantum Physics and Noncommutative Geometry}
\author{Lucian M. Ionescu}
\address{Mathematics Department\\Kansas State University\\
             Manhattan, Kansas 66502}
\email{luciani@@math.ksu.edu}
\date{v.1:06/99, v.2:06/00}

\begin{abstract}
The {\em quantum-event / prime ideal in a category/ noncommutative-point}
alternative to
{\em classical-event / commutative prime ideal/ point} is
suggested.

Ideals in additive categories, prime spectra and representation of quivers
are considered as mathematical tools appropriate to model quantum mechanics.
The space-time framework is to be reconstructed from the spectrum of
the path category of a quiver.

The interference experiment is considered as an example.
\end{abstract}

\maketitle
\tableofcontents


\section{Introduction}\label{S:intro}
Topological quantum field theory \cite{W1,At1,C,CF,CY1} is a
rigorous mathematical framework for the Feynman path integral formalism,
expressed in categorical language.
It is a ``categorification'' of quantum physics,
in the broader sense of the term \cite{CY2,Baez,I1}.
Together with the progress in quantum groups and geometric quantization,
there is enough evidence to say that, at a deeper level of understanding,
what is referred to in the physical realm as quantum phenomena,
corresponds in mathematics to noncommutative geometry...

... and, as stated in \cite{Ro1}, the viewpoint originally due to
Grothendieck and expressed in \cite{Ma} (p.83) is that
{\em ``To do geometry you really don't need a space.
All you need is a category of sheaves on this would-be space.''}

Our intention in this essay, is to bring together some ideas
and categorical tools, and to state some questions and comments on
the categorification of the formalism of quantum mechanics
in relation to non-commutative geometry.

\subsection{Space. Of what kind?}
Since 1930's general relativity, the physical theories based on
geometric methods evolved by adding {\em internal degrees of freedom}
to the {\em external degrees of freedom}
(adding ``structure'' to geometric points)
due to the presence of a space-time model, in an attempt to account for
the other fundamental interactions.

In Kaluza-Klein theories the additional dimensions were still
{\em external dimensions}.
Electromagnetism as a $U(1)$-Yang-Mills theory
has one additional {\em internal dimension}, nonabelian YM still add
extra {\em internal dimensions} etc.
One may ``extrapolate'' and ``predict'' that quantum physics of
the fundamental interactions should evolve
towards a no-external-all-internal dimensions type of theory.
More precisely, we think of quantum field theory being modeled
as an abstract theory from which space-time,
not anymore a fundamental concept,
is reconstructed not necessarily as a limit, but rather as a moduli space
from a categorical structure, as part of a {\em reduction principle}.

In this abstract theory the basic concept is not the
{\em space-time event} (classical event)
as a potential occurrence or ``existence'',
but rather an {\em elementary process} (quantum event).
The {\em interaction} itself should be the primary concept,
and not built on top of a preexisting structure.
In mathematical terms this should correspond to the natural
evolution from a ``set-theoretical approach'' 
with emphasis on ``elements + structure'' to a
categorical approach with its emphasis on {\bf morphisms}
and with the structure given through the category itself.

At the physical level, one may identify this tendency as part of the
process of deemphasizing configuration spaces,
and studying Poisson manifolds rather then cotangent spaces.
More then that, one should perhaps {\em model spaces of evolutions}
in a direct manner, rather then to represent them based on a state space
(e.g. as integral curves).

\subsection{Gravity. Too classical?}
The incompatibility between general relativity and quantum mechanics
is well-known and accounted for by the 50-60 years of
unconvincing attempts to ``quantize gravity'' (\cite{Po,BI}).
There must be a deep reason for the incompatibility
between general relativity, a classical-phenomenological-macroscopical
theory, and quantum physics, a fundamental-highly experimental theory,
and for the unsuccessful attempts to unify gravity with the other
fundamental interactions.
Perhaps general relativity is an effective theory,
as for example Fermi theory,
which is not 
``...cloaked in the beautiful geometric structure of General Relativity,
that may lead us astray in our search for a solution to the problem
of nonrenormalizability...'' \cite{Po}, p.74.

We would like to question our phenomenological understanding of space-time,
and recall that there are quantum gravity programs which reject
``the manifold conception of space-time'' \cite{BI}, p.4.
(see also \cite{Ng}, p.2).

We consider that the classical algebraic-geometry-physics point of view:
$$ commutative\ prime\ ideal \leftrightarrow \quad 
space-time\ point \quad\leftrightarrow \quad
classical\ event$$
should be replaced by:
$$ prime\ ideal\ in\ a\ category \quad\leftrightarrow\quad
NC-point\quad \leftrightarrow\quad
quantum \ event$$
An abstract C(ategorical)QFT,
based on representations of quivers or abstract categories,
and abandoning the manifold approach to space-time (compare \cite{BI}),
should generalize TQFTs based on a cobordism category
(see also \cite{Fu1}, Ex.2.2 and Remark 2.3, p.103).
The ``space-time'' should be reconstructed as a prime spectrum
of a category, as part of a ``correspondence principle'',
and should not be considered a primary concept.

As an ``example'' to give substance to this idea,
one may consider the Feynman diagrams as representations of quivers
bounding a certain basic quiver specific to the quantum experiment
being modeled, and determined by the ``geometry'' and
preparation stages of the experiment.
The allowable Feynman diagrams bounded by the basic quiver should be
related to the cut-off procedures applied to the experimental data.

The divergency of the Feynman integrals defined on the reconstructed
space, and due to the {\em unnecessary abundance of topological paths},
should be reduced according to a ``Hodge-principle''
(existence of a ``canonical representative'', e.g. harmonic form)
applied to homotopy classes of paths.

The physical counterpart would be the ``particle-wave duality'',
where the quantum process (say involving an electron)
behaves ``like a particle'' in an ``irreducible subprocess'',
means that the calculation can be classical, on a geometric path
(the ``Hodge representative''), while the overall amplitude
is a result of the interference due to the nontrivial homotopy structure
of the basic quiver (integration on the homotopy classes),
so that overall the quantum process behaves ``like a wave'',
i.e. exhibits interference.

Due to the believe that quantum phenomena should be modeled by
noncommutative geometry, we consider as a long-term goal to develop a
``finer'' axiomatic theory of quantum mechanics based on
spectra of categories of representations of quivers,
rather then based on Hilbert spaces and operators,
in a natural correspondence with noncommutative geometry
(categorical analog of string theory?).

A unification of the fundamental interactions
(without gravity if understood only by macroscopic evidence),
is then perhaps conceivable not by focusing on finding the
``overall gauge group'', but at a categorical level, 
and obtaining ``gravity'' by reduction at the classical level
(on the ``moduli space''), as a byproduct of the unification.

\section{NC geometry}\label{S:ncg}
Assuming that the ''categorification'' and ''Spec''
functors should commute, we translate the definition
of an affine scheme.


\subsection{Ideals and the prime spectrum}
Ideals in additive categories are discussed in \cite{I2}, and
in the references cited therein.

Let $\C{A}$ be a additive category, i.e. a category such that
for any two objects $A, B$ in $\C{A}$, $Hom(A,B)$
is an object in the category $\C{A}b$ of abelian groups.

A (bilateral) {\em ideal} in $\C{A}$  is a class of morphisms
(subfunctor, \cite{I2}, p.2), stable under addition and composition
whenever defined.
More precisely, a left ideal $I$ in $\C{A}$, is a family of abelian subgroups
$I(A,X)$ of $Hom(A,X)$ indexed by $A, X\in \C{A}$,
and such that the family of groups is stable under left composition with
arbitrary morphisms $\phi:X\to Y$:
$$\diagram
I(A,\phi)=Hom(A,\phi)_{|I(A,X)}                   &             &
\C{A} \ar@/^1pc/[r]^{I(A)} \ar@/_1pc/[r]_{I(B)} \ar@{}[r]|{\Downarrow I(f)} & \C{A}b\\
A \ar@{=>}[dr]_{I(A,Y)} \ar@{=>}[r]^{I(A,X)} & X \dto^{\phi} & I(B,X) \dto_{\phi\circ \cdot} \rto^{\cdot\circ f} & I(A,X) \dto^{\phi\circ\cdot} \\
 					     & Y             & I(B,Y) \rto_{\cdot\circ f}& I(A,Y)
\enddiagram$$
The double arrow denotes a set of morphisms, and the left diagram
is not necessarily ``commutative''.
Composing $\phi$ with morphisms from $I(A,X)$
may result in a proper subset of $I(A,Y)$.

The approach in \cite{Ro2} to define the spectrum of an abelian category
leads to a reconstruction theorem extending a know result of
Gabriel.
Since we consider that the primary concepts are the morphisms
we will adopt another approach.
We will consider the categorification of the usual definition
of a prime ideal \cite{AM}, \cite{Ja} (p.2).
For alternative definitions see \cite{Ba} (p.191),
\cite{SD} (p.140), \cite{BB} (p.1173).

An ideal $P$ is a {\em prime ideal} iff for any compassable morphisms
$f$ and $g$, $f\circ g\in P$ implies either $f\in P$ or $g\in P$.
The {\em prime spectrum} of $\C{A}$ is the class of prime ideals
(or category, if ideals are viewed as subfunctors),
with the usual Zarisky topology
(as a Grothendieck topology, to be made precise elsewhere).


\subsection{Representations of quivers and algebras}
Between sets and categories, a {\em quiver} $Q$ (or {\em precategory}) is a
``small category'' without a composition of morphisms.
A {\em morphism of quivers} is a ``functor'' without the 
condition regarding the composition of morphisms,
since there is none yet defined.

The category of quivers $\C{Q}v$ relates to the (2-)category of
small categories $\C{C}at$ similar to how sets relates to monoids.
There is a forgetful functor $U: \C{C}at\to \C{Q}v$,
with a left adjoint $\C{P}$, giving the free category
generated by a given quiver:
$$Hom_{\C{C}at}(\C{P}(Q),\C{C})\equiv Hom_{Qv}(Q, U(\C{C}))$$
$\C{P}(Q)$ is also called the {\em path category} generated by $Q$.

Let $\C{C}$ be an arbitrary category and $k$ a ring, e.g. $\B{Z}$.
Consider the category $\C{A}=k\C{C}$ with the same objects as $\C{C}$ and
with morphisms the free $k$-modules generated by the morphisms of the
category $\C{C}$.
Then $\C{A}$ is a $k$-linear category, and in particular additive,
in an obvious way.

To an arbitrary $k$-linear small category (additive, if $k=\B{Z}$),
one may associate a $k$-algebra (algebra),
called its {\em matrix algebra} (see \cite{Mi}, p.33).
The matrix algebra of the free $k$-linear category associated to the
path category generated by a quiver $Q$ is the {\em path algebra} $kQ$.
It is a categorification of the group ring construction,
corresponding to the ``multi-object'' case.
Explicitly, on generators, the product of two paths is zero,
unless they produce a concatenated new path \cite{Ri} (p.146).

Note that any additive category $\C{A}$ is a quotient of a path category
$\C{P}(U\C{A})$ by the kernel (ideal in $ \C{P}(U\C{A})$) of the
path evaluation (by composition) functor:

Alternatively, the elements of $C_\bullet(\C{A})=\C{P}\circ U (\C{A})$
may be viewed as singular simplices, i.e. functors from
the standard simplex (semisimplicial category) $\Delta_n$ to $\C{A}$
(\cite{I2}, p.8). In this way morphisms may be thought of as
1-simplices, and the path algebra, as a $k$-module,
is just the (graded) free $k$-module generated by singular simplices
of (categorical) singular homology.

A {\em representation of quivers} $F:Q\to U(k-mod)$ is a
morphism of quivers valued in a category of modules.
It extends uniquely to a functor $F:\C{P}(Q)\to k-mod$.
If $k=\B{Z}$, then it is a $\C{P}(Q)$-module \cite{St} (p.140),
\cite{Ba} (p.171).
A morphism of representations of quivers is a restriction
(with respect to morphisms) of a natural transformation of functors.

Regarding the quantum mechanics formalism,
considering representations of quivers rather then vector spaces
represents a gain in the initial structure involved.


\subsection{Algebras and Auslender-Reiten quivers}
The idea of generalizing the ordinary ``commutative'' geometry
by studying non-commutative algebras goes back to von Neumann
and his contemporaries \cite{W2} (p.253).
The inconvenience is that such an algebra might not have
enough ideals (``global'' ideal, contrasted to categorical ideals).
It is a ``pointless geometry''.

We point out that there is a well known correspondence between the
one object, ''total'' algebra, point of view and the
multi-objects point of view, where additive categories are
just ''rings with several objects'' \cite{Mi}
(groupoids = multi-object groups, etc.).
For a more detailed discussion of this correspondence, see \cite{Ri} and
\cite{AP,MP,Bl}.

We only note that to a finite dimensional algebra $R$
over an algebraic closed field $k$, one may associate
the Alexander-Reiten quiver, having vertices (objects)
isomorphism classes of indecomposable $R$-modules (i.e. $Spec(R-mod)$),
and arrows bases in a certain group of irreducible maps \cite{Ri} (p.181).
If $R$ is a basic algebra \cite{Ri} (p.147), then $R$
is of the form $k\Gamma/I$ with $\Gamma$ a quiver and $I$ an ideal
in $k\Gamma$.



\subsection{Points or paths?}
Feynman path integral formalism leads to mathematical difficulties,
one may argue, due to the philosophy ``{\bf space determines the paths}''.
For example, once a typical ``space-time '' is given, the set of
``acceptable paths'' is determined (topological, or smooth,
piecewise smooth, etc.), leading to an unnecessary large space
for integration in the lagrangean formalism.

Consistent with the emphasis of quantum processes / NC-geometry on:
$$ quantum\ event <-> interaction / correlation <-> path / morphism$$
it is legitimate to require having the ``points'' {\bf and}
``paths'' given at the same time, as a category for example.
The structure of the category encoded in morphisms and their composition
$\circ$, has a causal character (''local time''),
while a monoidal structure $\otimes$ models the ''spatial structure''.

\section{Quantum mechanics: an example}\label{S:qme}
\subsection{Interference experiment}
We will consider the elementary interference experiment with electrons
as explained in \cite{F} (p.37-5):

\epsfbox{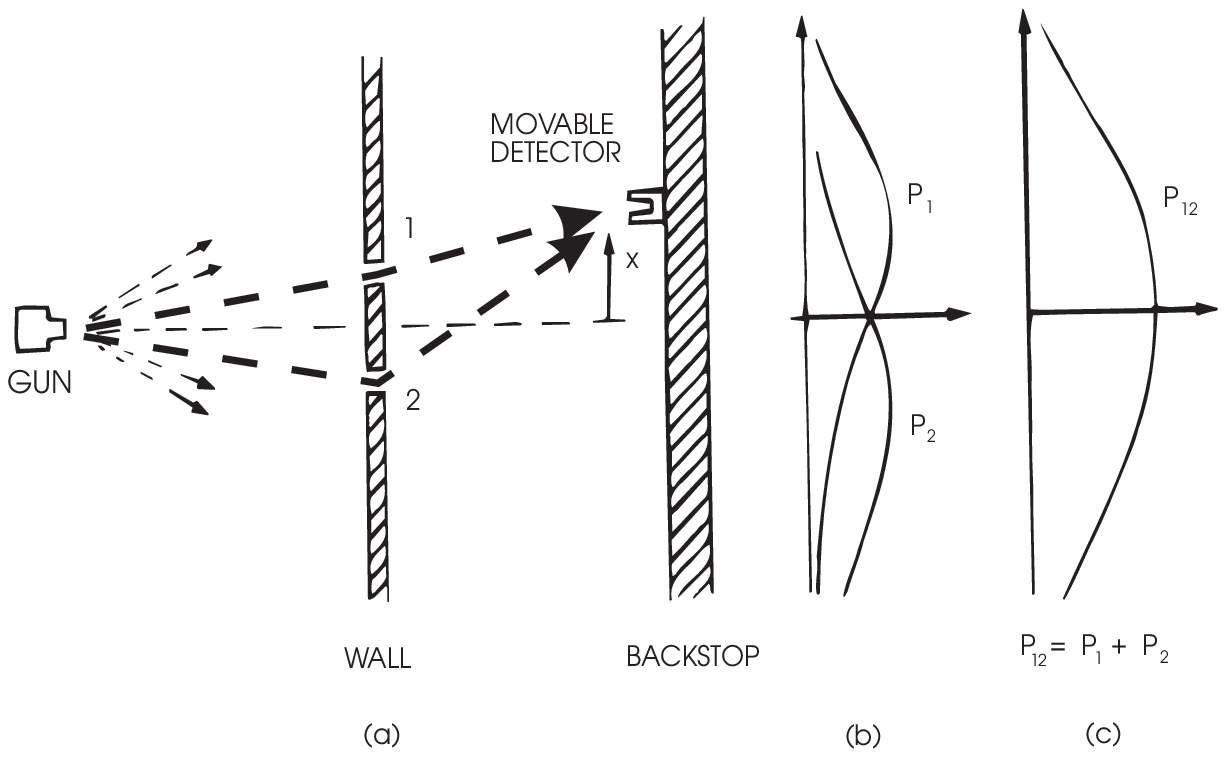}

from which the ``first principles of quantum mechanics'' are derived.

{\em ``An experiment is one in which all of the initial and
final conditions of the experiment are completely specified.''
What we will call {\bf an event} is, in general, just a specific set
of initial and final conditions.}
(\cite{F} p.37-10)

A categorical mind might choose to interpret initial and final conditions
as objects, more specifically sources and targets
(e.g. TQFTs, \cite{W1,At1,At2}).
We will consider the following principle:

To any {\em quantum experiment} there is an
{\em associated quiver} corresponding to the ``geometry'' of the experiment,
including the source(s), preparation stages and finally,
the measurement devices.

A perturbative approach would consider quivers (diagrams) containing it,
and should be thought of as a {\em bounding quiver}.
A filtering of the experimental data (cut-off procedure),
would alter the class of allowable quivers corresponding to
this bounding quiver.

For the above experiment, it is reasonable to prescribe the following
quiver $Q$ and representation $S$:
{\Large
$$\xymatrix @C=5pc @R=.1pc {
 & \overset{1}{\bullet} \drto^{} &\\
\overset{S}{\bullet} \urto^{} \drto^{} & & \overset{T}{\bullet}
 & \ar@/^/[r]^{{\bf S}} & & \B{C}\C{G}\\
& \overset{2}{\bullet}  \urto_{} &
} $$
}
A gauge group $G=U(1)$ interpreted as a one object category $\C{G}=\C{T}_G$,
and an action functor $\C{S}$ as a representation of the quiver
(valued in the additive category $\B{C}\C{G}$),
are also associated:
$$\C{S}(A)=*_G, \qquad \C{S}(\gamma)=e^{i \C{L}(\gamma)}$$
with $\C{L}$ a geometric lagrangean $\C{L}(\gamma)=length(\gamma)$.

{\em When an event can occur in several alternative ways,
the probability amplitude for the event is the sum of the probabilities
amplitudes for each way considered separately.}
(\cite{F} p.37-10)

The alternative ways for an  interaction, or ``paths'',
should be viewed as morphisms (e.g. cobordisms in TQFTs \cite{At1,At2}),
and the amplitude is:
$$ \int\limits_{Hom_\C{P}(S, T)} \C{S}(\gamma) \C{D}\gamma \quad = \quad
\int\limits_{Spec(\C{P})} \C{S}(P) \C{D}P
$$
where $\C{P}=\B{C}\C{P}(Q)$ is the $\B{C}$-linear category of paths
generated by $Q$. Here the integration is just the summation over the
``possible ways''  (morphisms) constituting the ``event''.

Alternatively, the amplitude can be viewed as an integral over
the prime spectrum of the path category,
and consisting of the two maximal ideals:
$$P_1=\B{C}\{S\overset{S1}{\to} 1,\ 1\overset{1T}{\to}  T,\ 1T\circ S1\}$$
$$P_2=\B{C}\{S\overset{S2}{\to} 2,\ 2\overset{2T}{\to}  T,\ 2T\circ S2\}$$
with the appropriate definition of the amplitude on each ``NC-point''.
This latter may be thought of as the result of enriching
the category $\C{P}$ in a category of Lagrange spaces,
or in a category of complexes in a way Fukaya categories are defined
(see \cite{Fu1} and references therein).

\subsection{Quantum versus cosmological event}
Although stated before, we would like to stress the following two levels.
One is the classical / macroscopical / phenomenological level:

\begin{center}

\vspace{.1in}
EVENT $\leftrightarrow$ Space-Time POINT
$\leftrightarrow$ (commutative ring) prime IDEAL
\end{center}

\vspace{.1in}
and the other is modern / quantum level / fundamental level:

\begin{center}
Quantum INTERACTION $\leftrightarrow$ NC POINT
$\leftrightarrow$ (categorical) prime IDEAL
\end{center}

\subsection{Particle-wave duality}\label{S:pwdp}
The modeling of the interference experiment as a quantum process and
studying the ``electron''
(or rather the phenomenon ``electron in this experimental setup''),
was done in section \S \ref{S:qme}
by associating the basic quiver and its representation (action).

It can be interpreted as considering the four {\em elementary processes}
on which the ``electron'' behaves like a particle,
and consequently defining the action $S(\gamma)$
as corresponding to a geometric path,
as if operating a ``homotopical'' reduction of the
space of all Feynman paths.
The overall behavior of the ``electron'',
i.e. the result of the interference due to a non-trivial $\pi_1$,
is the ``behavior like a wave'', within this experiment.
In this way the ``particle-wave duality'' principle is a natural
consequence of the theory.

\section{Conclusions}
Whether it is ``cohomological physics'' or ``categorical physics'',
it benefitted from more and more abstract mathematical models.

In this essay we have tried to raise some implicit questions,
through the introduction of some ideas and goals:

1) Formulation of quantum physics in categorical language,
in terms of representations of quivers.

2) The mathematical model should be relevant to
non-commutative geometry.

3) The primary concepts are Feynman diagrams, viewed as representations
of quivers.
Paths (morphisms) should be considered part of the primary  data,
along with kernel operators attached as representations of
quivers.

The reduction to homotopy classes should be viewed as a principle
establishing the correspondence with the traditional approach
to Feynman integrals.

4) Space-time is no longer a basic concept in quantum physics,
and should be reconstructed as a moduli space, rather then
as a ``classical limit''.

5) ``Particle-wave'' duality principle should be a natural consequence of
the theory. 
The ``wave-like behavior'' (interference) occurs as a consequence of the
multitude of ``homotopy classes'' of paths.
The ''particle-like behavior'' is the result of modeling
a subprocess as a classical propagation
along a classical (minimal/ geometrical) path.

6) Gravity should be obtained at a classical level as a byproduct of a
unification of the three fundamental interactions not by
enlarging the gauge group, but by enriching the categorical structure
of the theory.

To conclude,
as an alternative to the manifold approach to space-time,
representations of quivers should replace the classical
functions on configuration spaces, and pathes being thought of as
encoding the possible transitions of the system.
The ``reality'' should be less thought of as a classical-limit,
where one may forget its ``true quantum nature'',
and evaluations of partition and n-point functions
should rather be considered as ``quantum computing'',
where the mathematical quantum model is used both as a
theoretical computational device
to be ``built'' into a specific quantum computing device,
as well as a theory which should be tested on such a quantum computer
\cite{DE}.



\end{document}